\documentclass[a4paper,11pt]{article}

\usepackage{a4wide}
\usepackage{amsthm,amsmath,amssymb}

\def\C{\mathbb C}
\def\eps{\varepsilon}
\def\ga{\gamma}
\def\gb{{\gamma}^{-1}}
\def\gan{\gamma^{n}}
\def\gbn{{\gamma}^{-n}}
\def\GL{{\rm GL}}
\def\GLF{{\rm GL}_2(F)}

\def\I2{{\rm I}_2}
\def\In{{I_n}}
\def\ext{{\rm \bf ext}}

\def\MOF{{\rm M}_2({\mathcal O}_F)}
\def\N{\mathbb N}
\def\OF{{\mathcal O}_F}
\def\PF{\varpi}
\def\P{{\mathbb P}^1(F)}
\def\Q{\mathbb Q}
\def\sp{\rm Sp}
\def\res{{\rm \bf res}}

\def\v{{\rm val}}

\newtheorem{theo}{Theorem}

\newtheorem{lemma}{Lemma}

\begin{document}

\title{Test vectors for trilinear forms, given  two unramified representations}
\author{Louise Nyssen \cr {\footnotesize  lnyssen@math.univ-montp2.fr} }
\date{\today}

\maketitle

\section{Introduction}

\bigskip
Let $F$ be a finite extension of  $\Q_p$, with ring of integers $\OF$, and uniformizing parameter $\PF$, whose residual field has  
$q$ elements. For $G=\GLF$, let  $(\pi_1, V_1)$, $(\pi_2, V_2)$ and $(\pi_3, V_3)$ be three irreducible, admissible, infinite dimensional 
representations of  $G$. 
Using the theory of Gelfand pairs, Dipendra Prasad proves in \cite{P} that the space of $G$-invariant linear forms on $V_1\otimes V_2 \otimes V_3$ 
has dimension at most one. He gives a precise criterion for this dimension to be  one, that we will explain now.

Let $D_F^*$ be the group of invertible elements of  the quaternion division algebra  $D_F$ over $F$. 
When $(\pi_i, V_i)$ is a discrete series representation of $G$, denote by $(\pi'_i, V'_i)$  the irreducible representation of $D^*_F$ associated 
to $(\pi_i, V_i)$ by the Jacquet-Langlands correspondence. Again, by the theory of Gelfand pairs, 
the space of $D^*_F$-invariant linear forms on $V'_1\otimes V'_2 \otimes V'_3$ has dimension at most one.

Let $\sigma_i$ be the two dimensional representations of the Weil-Deligne group of $F$ associated to the irreducible representations $\pi_i$. 
The triple tensor product $\sigma_1 \otimes \sigma_2 \otimes\sigma_3$ is an eight dimensional symplectic representation of the Weil-Deligne group, 
and has local root number $\eps(\sigma_1 \otimes \sigma_2 \otimes \sigma_3)=\pm 1$. 
When $\eps(\sigma_1 \otimes \sigma_2 \otimes \sigma_3)=- 1$, 
one can prove that the representations $\pi_i$'s are all discrete series representations of $G$.

\begin{theo} (Prasad, Theorem 1.4 of \cite{P} ) Let   $(\pi_1, V_1)$, $(\pi_2, V_2)$, $(\pi_3, V_3)$ be three irreducible, admissible, infinite dimensional 
representations of  $G$ such that the product of their central characters is trivial. 
If all the representations $V_i$'s are cuspidal, 
assume that the residue characteristic of $F$ is not 2. Then \par
$\centerdot$ $\eps(\sigma_1 \otimes \sigma_2 \otimes \sigma_3)=1$ if and only if there exists a non-zero $G$-invariant linear form on 
$V_1\otimes V_2 \otimes V_3$ \par
$\centerdot$ $\eps(\sigma_1 \otimes \sigma_2 \otimes \sigma_3)=-1$ if and only if  there exists a non-zero $D_k^*$ invariant linear form on 
$V'_1\otimes V'_2 \otimes V'_3$.\par
\end{theo}

Given a non zero $G$-invariant linear form $\ell$ on $V_1\otimes V_2 \otimes V_3$, or a non-zero  $D_k^*$-invariant linear form $\ell'$ on 
$V'_1\otimes V'_2 \otimes V'_3$, the goal is to find a vector in  $V_1 \otimes V_2 \otimes V_3$ which is not in the kernel of $\ell$, 
or a vector in  $V'_1\otimes V'_2 \otimes V'_3$ which is not in the kernel of $\ell'$.  Such a vector is called a test vector. 
At first sight, it appears to have strong connections with the new vectors $v_1$, $v_2$ and $v_3$   of the representations $\pi_1$, $\pi_2$ and $\pi_3$.

\begin{theo}\label{vt-000} (Prasad, Theorem 1.3 of  \cite{P}) 
When  all the $\pi_i$'s are unramified principal series representations of $G$, $v_1 \otimes v_2 \otimes v_3$ is a test vector for $\ell$. 
\end{theo}

\begin{theo}\label{vt-111} (Gross and Prasad, Proposition 6.3 of \cite{GP}) When all   the $\pi_i$'s are  unramified twists of the special representation of $G$ : \par  

$\centerdot$ if $\eps(\sigma_1 \otimes \sigma_2 \otimes \sigma_3)=1$, then $v_1 \otimes v_2 \otimes v_3$ is a test vector for $\ell$,

$\centerdot$ if $\eps(\sigma_1 \otimes \sigma_2 \otimes \sigma_3)=-1$, let $R'$ be the unique maximal order in $D_F$. Then the open 
compact subgroup $R'^* \times R'^*\times R'^*$ fixes a unique line in  $V'_1\otimes V'_2 \otimes V'_3$. Any vector on this line is a test vector for $\ell'$. 
\end{theo} 

The proof by  Gross and Prasad of the first statement of this theorem actually contains another result:

\begin{theo}\label{vt-110} 
When two of the $\pi_i$'s are  unramified twists of the special representation of $G$ and the third
one belongs to the unramified principal series of $G$, $v_1 \otimes v_2 \otimes v_3$ is a test vector for $\ell$. 
\end{theo}

But the paper \cite{GP} gives evidence that $v_1 \otimes v_2 \otimes v_3$ is not always a test vector for $\ell$.
Let $K=\GL (\OF)$ be the maximal compact subgroup of $G$. If  $\pi_1$ and $\pi_2$ are unramified and if $\pi_3$ has conductor $n \geq 1$,
$\ell$ being $G$-invariant, $v_1$ and $v_2$ being $K$-invariant, one gets a $K$-invariant  linear form
$$\left\{ \begin{matrix}
V_3 & \longrightarrow & \C \hfill \cr
v & \longmapsto & \ell(v_1\otimes v_2\otimes v) \cr
\end{matrix}  \right.$$
which must be 0 since $\pi_3$ is ramified. Then $\ell(v_1\otimes v_2\otimes v_3)=0.$

Now Gross and Prasad make the following suggestion. Let $\In$ be the congruence subgroup
$$\In = \Bigl\{  \begin{pmatrix}
a & b \cr
c & d \cr
\end{pmatrix} \in K \quad  \vert \quad c \equiv 0 \quad mod \, {\PF}^n \quad \Bigr\}$$
and $R$ be a maximal order ${\rm M}_2(F)$ such that $R^* \cap K = \In$. If $v_1^*$ is a $R^*$-invariant vector in $V_1$, the linear form
$$\left\{ \begin{matrix}
V_3 & \longrightarrow & \C \hfill \cr
v & \longmapsto & \ell(v_1^*\otimes v_2 \otimes v) \cr
\end{matrix}  \right.$$
is invariant under the action of $R^* \cap K = \In$, and one can still hope that $v_1^* \otimes v_2 \otimes v_3$ is  a test vector for $\ell$.\par

\bigskip
The purpose of this paper is to prove that 
$v_1^* \otimes v_2 \otimes v_3$ actually {\it is}  a test vector for $\ell$. 
This is  the object of Theorem \ref{vt-00n}. 
The case $n=1$, together with Theorems \ref{vt-000}, \ref{vt-111} and \ref{vt-110}, 
 complete the study of test vectors when the $\pi_i$'s have ramification 0 or 1.\par

\bigskip
In the long term, the search for test vectors is motivated by the subconvexity problem for $L$-functions.
Roughly speaking, one wants to bound some $L$-functions along the critical line $\Re(z)=\frac{1}{2}$. 
A recent and successful idea in this direction has been to relate triple products of automorphic forms 
to special values of $L$-functions on the critical line. In \cite{BR1} and \cite{BR2} Joseph Bernstein and Andre Reznikov 
established a so called {\sl subconvexity bound} the  for the $L$-function of a triple of representations : 
each representation is attached to the eigenvalue of a certain operator, and  
the eigenvalue of one representation varies. Philippe Michel and Akshay Venkatesh 
considered the case when the {\sl level} of one representation varies. 
More details about subconvexity and those related techniques can be found in \cite{V} or \cite{MV}. 
Test vectors are key ingredients. Bernstein and Reznikov use an explicit test vector. 
Venkatesh uses a theoretical one, but explains that the bounds would be better with an explicit one (see paragraph 5 of \cite{V}). 
Unfortunately, the difficulty of finding them increases with the ramification of the representations involved.  

There is an extension of Prasad's result in \cite{HS}, where Harris and Scholl prove that 
the dimension of the space of $G$-invariant linear forms on $V_1\otimes V_2 \otimes V_3$ is one when $\pi_1$, $\pi_2$ and $\pi_3$ 
are principal series representations, either irreducible or reducible with their unique irreducible subspace, infinite dimensional. 
They apply the global setting of this to the construction of elements in the motivic cohomology of the product of two modular curves constructed by Beilinson. 

\bigskip
I would like to  thank Philippe Michel for suggesting this problem, 
Wen-Ching Winnie Li who invited me to spend one semester at PennState University where
I wrote the first draft of this paper, 
and of course Benedict Gross and Dipendra Prasad for the inspiration.  
I would also like to thank Paul Broussous and Nicolas Templier for many interesting discussions, 
and Eric Bahuaud for his help with English. 

 In a previous version of this paper, I obtained Theorem \ref{vt-00n} 
under an unpleasant technical condition. I am extremely grateful to Malden Dimitrov, 
because, thanks to our discussions on the subject, I found the way to remove the condition. 
In \cite{DN}, we are working on a more general version of  Theorem \ref{vt-00n}.

\section{Statement of the result}

\subsection{About induced and contragredient representations}\label{notations}

Let $(\rho, W)$ be a smooth representation of a closed subgroup $H$ of $G$. Let $\Delta_H$ be the modular function on $H$. 
The induction of $\rho$ from $H$ to $G$ is a representation $\pi$ whose space is the space 
${\rm Ind}_{H}^{G}\bigl( \rho \bigr)$  of functions $f$ from $G$ to $W$ 
satisfying the two following conditions :

(1) $\forall h \in H, \quad \forall g \in G, \quad f(hg)={\Delta_H}^{-\frac{1}{2}}(h) \rho(h) f(g)$,

(2) there exists an open compact subgroup $K_f$ of $G$ such that $$\forall k \in K_f, \quad \forall g \in G, \quad f(gk)= f(g)$$

\noindent where  $G$ acts by right translation. The resulting function will be denoted $\langle \pi(g) , f\rangle$ that is
$$\forall g, g_0 \in G, \quad \langle \pi(g),f\rangle(g_0) = f(g_0g).$$ 
With the  additional condition  that $f$ must be compactly supported modulo $H$, 
one gets the {\it compact} induction denoted by ${\rm ind}_{H}^{G}$. 
When $G/H$ is compact, there is no difference between ${\rm Ind}_{H}^{G}$ and ${\rm ind}_{H}^{G}$. 

Let $B$ the Borel subgroup of upper triangular matrices in $G$, and let $T$ be the diagonal torus. 
The character  $\Delta_T$ is trivial and we will use $\delta = {\Delta_B}^{-1}$ with 
$\delta \Bigl( \begin{pmatrix} a & b \cr 0 & d \cr \end{pmatrix} \Bigl)  = \vert \frac{a}{d} \vert$ 
where $ \vert \cdot \vert$ is the normalised valuation of $F$. The quotient $B \backslash G$ is compact and can be identified with $\P$.

For a smooth representation $V$ of $G$, $V^*$ is the space of linear forms on $V$. 
The contragredient representation $\widetilde{\pi}$ is given by the action of $G$ on  $\widetilde{V}$, the subspace of smooth vectors in $V^*$. 
If $H$ is a subgroup of $G$, $\widetilde{V} \subset \widetilde{V_{\vert H}} \subset V^*$. 

We refer the reader to  \cite{BZ} for more details about induced and contragredient representations.

\bigskip
\subsection{New vectors and ramification}\label{nv}
Let $(\pi,V)$ be an irreducible, admissible, infinite dimensional representation of  $G$ with central character $\omega$. 
To  the descending chain of  compact subgroups of $G$ 
$$ K \supset I_1 \supset \cdots \supset I_n \supset I_{n+1} \cdots  $$
one can associate an ascending chain of vector spaces
$$V^0=V^K \qquad {\rm and} \qquad \forall n\geq 1, \quad 
V^n = \Bigl\{  v \in V \quad \vert \quad   \forall  \begin{pmatrix} a & b \cr c & d \cr \end{pmatrix} \in I_n, \quad 
\pi \bigl( \begin{pmatrix} a & b \cr c & d \cr \end{pmatrix}\bigr) v = \omega(a) v 
\quad   \Bigr\}. $$
There exists a minimal $n$ such that the vector space $V^n$ is not $\{0\}$. 
It is necessarily one dimensional and any generator of $V^n$ is called a new vector of $(\pi, V)$. 
The integer $n$ is the conductor of $(\pi, V)$. 
The representation $(\pi,V)$ is said to be unramified when $n=0$. Else, it is ramified. 

More information about new vectors can be found in \cite{C}.

\bigskip
\subsection{The main result}\label{result}
Let  $(\pi_1, V_1)$, $(\pi_2, V_2)$ and $(\pi_3, V_3)$ be three irreducible, admissible, infinite dimensional 
representations of  $G$ such that the product of their central characters is trivial. 
Assume that $\pi_1$ and $\pi_2$ are unramified principal series, and that $\pi_3$ has conductor $n\geq 1$. 
According to Theorem 1, since $\pi_1$ and $\pi_2$ are not discrete series, 
there exists a non-zero, $G$-invariant linear form $\ell$ on $V_1\otimes V_2 \otimes V_3$. 
We are looking for a vector $v$ in $V_1\otimes V_2 \otimes V_3$ which is not in the kernel of $\ell$. 
In order to follow the suggestion of Gross and Prasad we consider

$$\ga=
\begin{pmatrix}
\PF & 0 \cr
0 & 1 \cr
\end{pmatrix} \qquad {\rm and} \qquad R_n = {\gbn} \MOF \gan .$$

One can easily check that 

$$R^*_n = \gbn K \gan \qquad {\rm and} \qquad R^*_n \cap K = \In .$$

If $v_1$, $v_2$ and $v_3$ denote the new vectors of $\pi_1$, $\pi_2$ and $\pi_3$, the vector
$$v_1^*=  \pi_1(\gbn) \cdot v_1 $$
is invariant under the action of  $R^*_n$. Hence we can write 

$$ v_1^*\in {V_1}^{R^*_n},  \qquad \qquad v_2 \in {V_2}^{K}  \qquad {\rm and} \qquad  v_3 \in {V_3}^{R^*_n \cap K }  .$$

\begin{theo}\label{vt-00n} Under those conditions, $v_1^*\otimes v_2 \otimes v_3$ is a test vector for $\ell$.
\end{theo}

The proof will follow the same pattern as Prasad's proof of Theorem \ref{vt-000} in \cite{P}, 
with the necessary changes.

\section{Going down Prasad's exact sequence}

\subsection{Central characters}\label{caractere-central}

Let $\omega_1$, $\omega_2$ and $\omega_3$ be the central characters of $\pi_1$, $\pi_2$ and $\pi_3$. 
Notice that the condition $\omega_1 \omega_2 \omega_3 =1 $ derives from the $G$-invariance of $\ell$.
Since $\pi_1$ and $\pi_2$ are unramified, $\omega_1$ and $\omega_2$ are unramified too, and so is $\omega_3$
because $ \omega_1 \omega_2 \omega_3 =1 $. 
Let $\eta_i$, for $i\in \{ 1,2,3 \}$ be unramified quasi-characters of $F^*$ with 
$\eta_i^2=\omega_i$ and $\eta_1 \eta_2 \eta _3 =1$. Then
$$V_1 \otimes V_2 \otimes V_3 \simeq 
\bigl(  V_1 \otimes \eta_1^{-1} \bigr) \otimes \bigl(V_2\otimes \eta_2^{-1} \bigr) \otimes \bigl(V_3\otimes \eta_3^{-1} \bigr) $$
as a representation of $G$. Hence it is enough to prove  Theorem \ref{vt-00n} when the central characters of the representations are trivial.

When $n=1$, it is also enough to prove Theorem \ref{vt-00n} when $V_3$ is the special representation $\sp$ of $G$ : 
take $\eta_3$ to be the unramified character such that $V_3 = \eta_3 \otimes \sp$.

\subsection{Prasad's exact sequence}\label{suites-exactes}

Let us now explain how Prasad finds $\ell$. It is equivalent to search for $\ell$ or to search for a non-zero element in   
${\rm Hom}_G \Bigl( V_1 \otimes V_2 , \widetilde{V_3} \Bigr)$. 
Since the central characters of $\pi_1$ and $\pi_2$ are trivial, there are unramified characters $\mu_1$ and $\mu_2$ such that  for $i=1$ and $i=2$
 $$\pi_i = {\rm Ind}_{B}^{G} \chi_i \qquad {\rm with} \qquad 
 \chi_i \Bigl( \begin{pmatrix} a & b \cr 0 & d \cr \end{pmatrix} \Bigl) = \mu_i \Bigl(\frac{a}{d} \Bigr).$$

\noindent Hence
$$V_1 \otimes V_2 = {\rm Res}_{G}\,{\rm Ind}_{B \times B}^{G \times G} \Bigl( \chi_1 \times \chi_2 \Bigr)$$
where $G$ is diagonally embedded in $G\times G$ for the restriction. 
The action of $G$ on $B\times B \backslash G\times G = \P \times \P$ has precisely two orbits.   
The first is  $\{ (u,v) \in \P \times  \P \quad \vert \quad u \not = v  \}$ which is open and  can be identified with $T \backslash G$.  
The second orbit  is  the diagonal embedding of $\P$ in $\P \times \P$, which is closed and can be identified with $B \backslash G$. 
Then, we have a short exact sequence of $G$-modules

\begin{equation}\label{courtesuite}
0 \rightarrow {\rm ind}_{T}^{G}\Bigl( \frac{\chi_1}{\chi_2} \Bigr)  \xrightarrow{\ext} V_1 \otimes V_2 \xrightarrow{\res} 
{\rm Ind}_{B}^{G} \Bigl(\chi_1 \chi_2 \delta^{\frac{1}{2}}  \Bigr)\rightarrow 0 .
\end{equation}

The surjection $\res$ is the restriction of functions from  $G \times G$ to the diagonal part of $B \backslash G \times B \backslash G$, that is
$$\Delta_{B \backslash G} = \Bigl\{ (g,bg) \quad \vert \quad b \in B, \quad g \in G \Bigr\}.$$ 
The injection $\ext$ takes a function 
$f \in {\rm ind}_{T}^{G}\Bigl( \frac{\chi_1}{\chi_2} \Bigr)$ to a function 
$F \in {\rm Ind}_{B \times B}^{G \times G} \Bigl( \chi_1 \times \chi_2 \Bigr)$ vanishing on $\Delta_{B \backslash G}$, 
and is given by the relation
$$F \Bigl( g, \begin{pmatrix} 0 & 1 \cr 1 & 0 \cr \end{pmatrix} g  \Bigr) = f(g) \label{rel},$$
on the other orbit. Applying the functor ${\rm Hom}_G \Bigl( \cdot \, , \widetilde{V_3} \Bigr) $, one gets a long exact sequence

\begin{multline}\label{longuesuite} 
0 \rightarrow {\rm Hom}_G \Bigl( {\rm Ind}_{B}^{G} \Bigl(\chi_1 \chi_2 \delta^{\frac{1}{2}}  \Bigr) , \widetilde{V_3} \Bigr) 
  \rightarrow  {\rm Hom}_G \Bigl( V_1 \otimes V_2 , \widetilde{V_3} \Bigr) 
  \rightarrow  {\rm Hom}_G \Bigl( {\rm ind}_{T}^{G}\Bigl( \frac{\chi_1}{\chi_2} \Bigr), \widetilde{V_3} \Bigr) \\
\hfill \downarrow \hskip20mm\\
\hfill \cdots  \leftarrow {\rm Ext}_G^1 \Bigl( {\rm Ind}_{B}^{G} \Bigl(\chi_1 \chi_2 \delta^{\frac{1}{2}}  \Bigr) , \widetilde{V_3} \Bigr)
\end{multline}

\subsection{The simple case}\label{cas-simple}

The situation is easier when $n=1$ and $\mu_1 \mu_2 \vert\cdot\vert^{ \frac{1}{2}} = \vert\cdot\vert^{- \frac{1}{2}}$, 
as $\pi_3$ is special and there is a natural surjection  
$${\rm Ind}_{B}^{G} \Bigl(\chi_1 \chi_2 \delta^{\frac{1}{2}}  \Bigr) \longrightarrow \widetilde{V_3}$$
whose kernel is the one dimensional subspace of constant functions. Thanks to the exact sequence (\ref{courtesuite}) one gets a surjection $\Psi$
$$\begin{matrix}
 V_1 \otimes V_2  & \xrightarrow{\res}  & {\rm Ind}_{B}^{G} \Bigl(\chi_1 \chi_2 \delta^{\frac{1}{2}}  \Bigr) \cr
\hfill {\scriptstyle \Psi}\! \searrow  & & \swarrow \hfill\cr
&\widetilde{V_3}&\cr
\end{matrix}$$
which corresponds to 
$$\ell\left\{ \begin{matrix}
 V_1 \otimes V_2 \otimes V_3 & \longrightarrow & \C \hfill \cr
 v \otimes v' \otimes v'' & \longmapsto & \Psi(v\otimes v') . v'' \cr
\end{matrix}\right.$$ 

The surjection $\Psi$ vanishes on $v_1^*\otimes v_2$ if and only if $\res(v_1^*\otimes v_2)$ has constant value on $\P \simeq B \backslash G$. 
An easy computation proves that $\res(v_1^*\otimes v_2)$ is not constant : 
the new vectors $v_1$ and $v_2$ are functions from $G$ to $\C$ such that  

$$\forall i \in \{1,2\}, \quad \forall b \in B, \quad \forall k \in K, \quad \quad v_i(bk)= \chi_i(b) \cdot \delta(b)^\frac{1}{2}$$  
and 
$$ \forall g \in G,  \quad \quad v_1^*(g)=v_1(g\gb). $$
Then
$$ ( v_1^* \otimes v_2 ) \Bigl( \begin{pmatrix}
1 & 0 \cr
0 & 1 \cr
\end{pmatrix} \Bigr)  
= v_1\Bigl( \gb \Bigr)  v_2\Bigl( \begin{pmatrix}
1 & 0 \cr
0 & 1 \cr
\end{pmatrix} \Bigr)
= v_1\Bigl(  
\begin{pmatrix}
\PF^{-1} & 0 \cr
0 & 1 \cr
\end{pmatrix}
\Bigr)  =  \mu_1(\PF)^{-1}\vert \PF \vert^{-\frac{1}{2}} =  \frac{\sqrt{q}}{\mu_1(\PF)} $$
and 
$$( v_1^* \otimes v_2 ) 
\Bigl( \begin{pmatrix}
0 & 1 \cr
1 & 0 \cr
\end{pmatrix} \Bigr)  
=  v_1\Bigl(  \begin{pmatrix}
0 & 1 \cr
1 & 0 \cr
\end{pmatrix}
\begin{pmatrix}
\PF^{-1} & 0 \cr
0 & 1 \cr
\end{pmatrix}
\Bigr)  
= v_1\Bigl(  \begin{pmatrix}
1 & 0 \cr
0 & \PF^{-1} \cr
\end{pmatrix}
\begin{pmatrix}
0 & 1 \cr
1 & 0 \cr
\end{pmatrix}
\Bigr)
=  \frac{\mu_1(\PF)}{\sqrt{q}}. $$

\noindent The representation $\pi_1$ is principal so $\frac{\sqrt{q}}{\mu_1(\PF)} \not = \frac{\mu_1(\PF)}{\sqrt{q}}$ and 
$$ ( v_1^* \otimes v_2 ) \Bigl( \begin{pmatrix}
1 & 0 \cr
0 & 1 \cr
\end{pmatrix} \Bigr)  \not=   ( v_1^* \otimes v_2 ) 
\Bigl( \begin{pmatrix}
0 & 1 \cr
1 & 0 \cr
\end{pmatrix} \Bigr).$$  
Hence $\Psi$ does not  vanish on $v_1^*\otimes v_2$. 
Now, $v_1^*$ being $R^*_1$-invariant and $v_2$ being $K$-invariant, $\Psi(v_1^* \otimes v_2)$ is a non-zero $I_1$-invariant element of $\widetilde{V_3}$, 
that is, a new vector for $\widetilde{\pi_3}$. Consequently it does not vanish on $v_3$ :
$$\ell(v_1^*\otimes v_2 \otimes v_3)= \Psi(v_1^* \otimes v_2).v_3 \not= 0$$
and $v_1^*\otimes v_2 \otimes v_3$ is a test vector for $\ell$.

\subsection{The other case}\label{autre-cas}

If $n \geq 2$ or $\mu_1 \mu_2 \vert\cdot\vert^{ \frac{1}{2}} \not = \vert\cdot\vert^{- \frac{1}{2}}$
then ${\rm Hom}_G \Bigl( {\rm Ind}_{B}^{G} \Bigl(\chi_1 \chi_2 \delta^{\frac{1}{2}}  \Bigr) , \widetilde{V_3} \Bigr) =0$ and by Corollary 5.9 of \cite{P} 
$${\rm Ext}_G^1 \Bigl( {\rm Ind}_{B}^{G} \Bigl(\chi_1 \chi_2 \delta^{\frac{1}{2}}  \Bigr) , \widetilde{V_3} \Bigr)=0.$$ 
Through the long exact sequence (\ref{longuesuite}) we get an isomorphism 
$${\rm Hom}_G \Bigl( V_1 \otimes V_2 , \widetilde{V_3} \Bigr) \simeq {\rm Hom}_G \Bigl( {\rm ind}_{T}^{G}\Bigl( \frac{\chi_1}{\chi_2} \Bigr), \widetilde{V_3} \Bigr),$$
and by Frobenius reciprocity 
$${\rm Hom}_G \Bigl( {\rm ind}_{T}^{G}\Bigl( \frac{\chi_1}{\chi_2} \Bigr) , \widetilde{V_3} \Bigr) \simeq 
{\rm Hom}_T \Bigl( \Bigl( \frac{\chi_1}{\chi_2} \Bigr) , \widetilde{{V_3}_{\vert T}} \Bigr),$$
where $\widetilde{{V_3}_{\vert T}}$ is the space of the contragredient representation of ${\pi_3}_{\vert T}$.
By Lemmas 8 and 9 of \cite{W}, the latter space is one dimensional. 
Thus, we have a chain of isomorphic one dimensional vector spaces
\smallskip
\begin{footnotesize}
$$\begin{matrix}
                 {\rm Hom}_G  \Bigl( V_1 \otimes V_2 \otimes V_3 , \C \Bigr) 
& \tilde{\rightarrow}  & {\rm Hom}_G \Bigl( V_1 \otimes V_2 , \widetilde{V_3} \Bigr) 
& \tilde{\rightarrow}  & {\rm Hom}_G \Bigl( {\rm ind}_{T}^{G}\Bigl( \frac{\chi_1}{\chi_2} \Bigr), \widetilde{V_3} \Bigr)
& \tilde{\rightarrow}  & {\rm Hom}_T \Bigl( \Bigl( \frac{\chi_1}{\chi_2} \Bigr) , \widetilde{{V_3}_{\vert T}} \Bigr)\cr
&&&&&&\cr
\ell & \mapsto & \Psi & \mapsto & \Phi & \mapsto & \varphi \cr
\end{matrix}\label{chaine}$$
\end{footnotesize}
with generators $\ell$, $\Psi$, $\Phi$ and $\varphi$ corresponding via the isomorphisms. 
Notice that $\varphi$ is a linear form on $V_3$ such that  
\begin{equation}\label{phi}
\forall t \in T, \qquad \forall v \in V_3, \qquad \varphi\bigl(\pi_3(t)v\bigr) = \frac{\chi_2(t)}{\chi_1(t)}\varphi(v) 
\end{equation}
which is identified to the following element of ${\rm Hom}_T \Bigl( \Bigl( \frac{\chi_1}{\chi_2} \Bigr) , \widetilde{{V_3}_{\vert T}} \Bigr)$
$$\left\{\begin{matrix} \C & \longrightarrow & \widetilde{{V_3}_{\vert T}}\cr
z & \longmapsto & z\varphi\cr
\end{matrix}\right.  $$

\begin{lemma}\label{lemmeGP}
$\varphi(v_3) \not= 0$.
\end{lemma}

{\it Proof} : this is Proposition 2.6 of \cite{GP} with the following translation : \par
- the local field $F$ is the same, \par
- the quadratic extension $K/F$ of Gross and Prasad is $F \times F$  
and their group $K^*$ is our torus $T$, \par
- the infinite dimensional representation $V_1$ of  Gross and Prasad is our $\pi_3$,\par
- the one dimensional, unramified representation $V_2$ of  Gross and Prasad is $\frac{\chi_1}{\chi_2}$.\par

Then the representation that Gross and Prasad call $V$ is $\frac{\chi_1}{\chi_2} \otimes \pi_3$ and their condition (1.3) is exactly our condition (\ref{phi}).
In order to apply Gross and Prasad's Proposition, we need to check the equality
$$\eps(\sigma \otimes\sigma_3 ) = \alpha_{K/F}(-1)\omega(-1).$$
Basically, it is true because $K$ is not a field. Let us give some details.  \par
- In \cite{GP}, $\omega$ is the central character of the representation $V_1$ which is trivial for us. \par
- The character $\alpha_{K/F}$ is the quadratic character of $F^*$ associated to the extension $K/F$ 
by local class-field theory. Here, it is trivial because $K$ is $F\times F$. \par
- To compute $\eps(\sigma \otimes\sigma_3 )$ we will use the first pages of \cite{T}. 
$$\forall \begin{pmatrix} x & 0 \cr 0 &  z\cr \end{pmatrix} \in T \qquad 
\frac{\chi_1}{\chi_2}\Bigl(\begin{pmatrix} x & 0 \cr 0 &  z \cr \end{pmatrix}\Bigr)
 = \frac{\mu_1}{\mu_2}(x)\frac{\mu_2}{\mu_1}(z) \Rightarrow 
\eps(\sigma \otimes\sigma_3 ) 
 = \eps(\frac{\mu_1}{\mu_2} \otimes\sigma_3 )\eps(\frac{\mu_2}{\mu_1} \otimes\sigma_3 )$$
Since the determinant of $\sigma_3$ is the central character of $\pi_3$ which is trivial, 
$\sigma_3$ is isomorphic to its own contragredient and the contragredient representation of 
$\frac{\mu_1}{\mu_2} \otimes\sigma_3$ is $\frac{\mu_2}{\mu_1} \otimes\sigma_3$. 
Formula (1.1.6) of \cite{T} leads to 
$$\eps(\sigma \otimes\sigma_3 )= {\rm det}\bigl(\sigma_3(-1)\bigr)=1= \alpha_{K/F}(-1)\omega(-1) .$$

\noindent According to \cite{GP}, the restriction of $\frac{\chi_1}{\chi_2} \otimes \pi_3$ to the group 
$$M=\Bigl\{  \begin{pmatrix} x & 0 \cr 0 & z \cr \end{pmatrix} \quad \vert  \quad x,y \in \OF^* \Bigr\} \times \In$$
fixes a unique line in $V_3$ : it is the line generated by the new vector $v_3$. 
Still according to  Gross and Prasad, a non-zero linear form on $V_3$ which satisfies (\ref{phi}) cannot vanish on $v_3$. 
$\hfill \Box$

\bigskip
\bigskip
We will deduce from lemma \ref{lemmeGP} that $\ell(v_1^*\otimes v_2 \otimes v_3) \not= 0$.

\section{Going up Prasad's exact sequence}

\subsection{From $\varphi(v_3)$ to $f$}\label{phi-f}

Let $f$  be the element of ${\rm ind}_{T}^{G}\Bigl( \frac{\chi_1}{\chi_2} \Bigr)$ 
which is the characteristic function of the orbit of the unit 
in the decomposition of $T \backslash G$ under the action of $\In$. This means~: 
\begin{equation}\label{f}
f(g) = \left\{ 
\begin{matrix}
\frac{\chi_1(t)}{\chi_2(t)} & {\rm if} \quad g=tk \quad {\rm with} \quad t \in T  \quad {\rm and} \quad k \in \In  \cr
0 & {\rm else} \hfill\cr
\end{matrix}
\right.
\end{equation}

Then, the function
$$\left\{ \begin{matrix}
 G & \longrightarrow & \C \hfill \cr
 g & \longmapsto &  f(g) \, \varphi\Bigl( \pi_3(g)v_3 \Bigr)\cr
\end{matrix}\right.   $$
is invariant by the action of $T$ by left translation and  we can do the following computation~: 
\begin{align*} 
\Bigl( \Phi(f) \Bigr)(v_3) & = \int_{T \backslash G} \! f(g) \, \varphi\Bigl( \pi_3(g)v_3 \Bigr) dg  \\
                            & =\int_{(T\cap K) \backslash\In} \! \varphi\Bigl( \pi_3(k)v_3 \Bigr) dk \\
                            & = \lambda \cdot \varphi(v_3) .\\
\end{align*}  
where $\lambda$ is a non-zero constant. Thanks to Lemma \ref{lemmeGP} we know that $\varphi(v_3) \not= 0$, so 
\begin{equation}\label{phi-f-v_3}
\Bigl( \Phi(f) \Bigr) (v_3) \not= 0. 
\end{equation}

\subsection{From $f$ to $F$}

Now, we are going to compute $F=\ext(f)$ in $V_1\otimes V_2$.
Let $a$ and $b$ be the numbers
$$a= \frac{\mu_1(\PF)}{\sqrt{q}} \qquad \qquad \qquad  b=\frac{\mu_2(\PF)}{\sqrt{q}}$$ 
They verify $$(a^2-1)(b^2-1) \not= 0$$
because $\pi_1$ and $\pi_2$ are  principal series representations. 
For the sake of simplicity, we shall use the following notation : for any $g$ in $G$ 
$$g v_1^* = \Bigl< \pi_1 (g), v_1^*  \Bigr> \qquad {\rm and} \qquad g v_2   =  \Bigl<  \pi_2 (g),  v_2 \Bigr>$$
  
\begin{lemma} \label{FV} The function $F$ is given by the formula 
$$F =   A \cdot v_1'\otimes v_2'$$
with
$$A = \frac{a^n}{(a^2-1)(b^2-1)}, \qquad 
v_1'= a \cdot \ga^{-(n-1)}v_1 - \ga^{-n}v_1 \qquad{\rm and}\qquad 
v_2'= b \cdot \gb v_2  -  v_2.$$
\end {lemma}

{\it Proof} :  the function  $f$ is described by formula (\ref{f}), 
and $\ext(f)$ is described by the short exact sequence (\ref{courtesuite}) using the orbits of the action of $G$ on $B\times B \backslash G \times G$.
The function $F$ must vanish on the closed orbit
$$\Delta_{B \backslash G} = \Bigl\{ (g,bg) \quad \vert \quad b \in B, \quad g \in G \Bigr\}.$$
The open orbit can be identified with $T \backslash G$ via the bijection 
$$\left\{ \begin{matrix}
T \backslash G & \longrightarrow & \Bigl( B \backslash G \times B \backslash G \Bigr) \setminus \Delta_{B \backslash G} \cr
\hfill Tg & \longmapsto & \Bigl( Bg , B\begin{pmatrix} 0 & 1 \cr 1 & 0 \cr \end{pmatrix}g \Bigr) \hfill\cr
\end{matrix}\right.$$
through which, the orbit of the unit in  $T \backslash G$ under the action of $\In$ corresponds to 
$$ \Bigl\{  \Bigl( Bk , B\begin{pmatrix} 0 & 1 \cr 1 & 0 \cr \end{pmatrix}k \Bigr) \quad \vert \quad k \in \In \Bigr\}  .$$
Now, pick any $(k,k') \in K \times K$. If $k' \in Bk$, then 
$$k \in \In \iff k' \in \In \qquad {\rm and} \qquad k \in I_1 \iff k' \in I_1.$$
When $k'\not\in Bk$, write 
$k=\begin{pmatrix}x & y \cr z & t\cr \end{pmatrix}$ and $k'=\begin{pmatrix}x' & y' \cr z' & t'\cr \end{pmatrix}$. 
There exists $(b_1, b_2) \in B \times B$ such that 
$$k=b_1k_0 \qquad{\rm and}\qquad k'=b_2 \begin{pmatrix} 0 & 1 \cr 1 & 0 \cr \end{pmatrix}k_0 \qquad {\rm with} \qquad 
k_0=\begin{pmatrix}z' & t' \cr z & t\cr \end{pmatrix} \in \MOF.$$
Then $(k,k')$ is in the orbit of the unit in  $T \backslash G$ under the action of $\In$ if and only if $k_0$ is in $T\In$. 
Because $k$ and $k'$ are in $K$, one can see that 
$$k_0 \in T\In  \quad  \iff \quad k_0 \in \In$$
and
\begin{align*}
k_0 \in \In  & \iff  z \equiv 0 \quad mod \, {\PF}^n \quad{\rm and}\quad z't \in \OF^* \cr
             &  \iff z \equiv 0 \quad mod \, {\PF}^n \quad{\rm and}\quad z' \in  \OF^*   \cr
             &   \iff k \in \In \quad{\rm and}\quad k' \notin I_1. \cr
\end{align*}

It follows that $(k,k')$ corresponds to an element of the orbit of the unit 
in the decomposition of $T \backslash G$ under the action of $\In$ if and only if 
$k \in \In$ and $k' \notin I_1$. Then, it will be enough to check that 
\begin{equation}\label{FK}
F(k,k') = \left\{ 
\begin{matrix}
1 & {\rm if} \quad k \in \In  \quad {\rm and} \quad  k' \not\in I_1 \cr
0 & {\rm else} \hfill\cr
\end{matrix}
\right.
\end{equation}
This is mere calculation. With $a= \frac{\mu_1(\PF)}{\sqrt{q}}$, we need

\begin{lemma}\label{calcul} 
Let $i$ be any element of $\N$. The values of the function $\ga^{-i}v_1$ on $K$ are given by the formula
\begin{equation}
\forall k =\begin{pmatrix}x & y \cr z & t\cr \end{pmatrix} \in K , \qquad \qquad 
\ga^{-i}v_1(k) = \left\{ 
\begin{matrix}
a^{i} \hfill & {\rm if} \quad \v(\frac{z}{t})\leq 0 \hfill \cr
& \cr
a^{i-2 \, \v(\frac{z}{t})} & {\rm if} \quad 1 \leq \v(\frac{z}{t})\leq i-1 \cr
& \cr
a^{-i} \hfill & {\rm if} \quad i \leq \v(\frac{z}{t}) \hfill \cr
\end{matrix}
\right.
\end{equation}

\end{lemma}

{\it Proof} :  first, recall that the new vector $v_1$is a function from $G$ to $\C$ such that  
$$\quad \forall b \in B, \quad \forall k \in K, \quad \quad v_1(bk)= \chi_1(b) \cdot \delta(b)^\frac{1}{2}.$$  
Then, for $k =\begin{pmatrix}x & y \cr z & t\cr \end{pmatrix}$ in $K$, 
either $z$ or $t$ is in $\OF^*$. The other one is in $\OF$. Write 
$$ k \ga^{-i}= \begin{pmatrix}
\frac{x}{\PF^i} & y \cr
\frac{z}{\PF^i} & t \cr
\end{pmatrix}.$$ 
If $\v\,\frac{z}{t}\leq 0$, then $\v\,z=0 $,  $\v(\frac{\PF^i t}{z})\geq 0$ and
$$ k \ga^{-i} =  \begin{pmatrix}
\frac{xt-yz}{z} & \frac{x}{\PF^i} \cr
0 &  \frac{z}{\PF^i} \cr
\end{pmatrix}\begin{pmatrix}
0 & -1 \cr
1 & \frac{\PF^i t}{z} \cr
\end{pmatrix}$$
with $$ \begin{pmatrix}
0 & -1 \cr
1 & \frac{\PF^i t}{z} \cr
\end{pmatrix} \in K \quad {\rm and} \quad (\chi_1 \cdot {\delta}^\frac{1}{2})\begin{pmatrix}
\frac{xt-yz}{z} & \frac{x}{\PF^i} \cr
0 & \frac{z}{\PF^i} \cr
\end{pmatrix} = \Bigl(\frac{\mu_1(\PF)}{\sqrt{q}}\Bigr)^{i-2\v\,z}=a^{i-2\v\,z}=a^i.$$ 
If $1\leq \v\,\frac{z}{t} \leq i-1$, then $\v\,t = 0 $, $\v(\frac{\PF^i t}{z})\geq 0$ and
the computation is quite the same, except that 
$$(\chi_1 \cdot {\delta}^\frac{1}{2})\begin{pmatrix}
\frac{xt-yz}{z} & \frac{x}{\PF^i} \cr
0 & \frac{z}{\PF^i} \cr
\end{pmatrix} = \Bigl(\frac{\mu_1(\PF)}{\sqrt{q}}\Bigr)^{i-2\v\,z}=a^{i-2\v\,\frac{z}{t}}.$$ 
If $i \leq \v\,\frac{z}{t}$, then $\v\,t = 0 $, $\v(\frac{\PF^i t}{z})\leq 0$ and
$$  k \ga^{-i}  =  \begin{pmatrix}
\frac{xt-yz}{t \,\PF^i } & y \cr
0 &  t \cr
\end{pmatrix}\begin{pmatrix}
1 & 0 \cr
\frac{z}{t \,\PF^i } &  1 \cr
\end{pmatrix}$$
with $$ \begin{pmatrix}
1 & 0 \cr
\frac{z}{t \,\PF^i } &  1 \cr
\end{pmatrix} \in K \quad {\rm and} \quad (\chi_1 \cdot {\delta}^\frac{1}{2})(\begin{pmatrix}
\frac{xt-yz}{t\, \PF^i } & y \cr
0 &  t \cr
\end{pmatrix}) = \Bigl(\frac{\mu_1(\PF)}{\sqrt{q}}\Bigr)^{-i-2\v\,t}=a^{-i}.$$ 
$\hfill\Box$

\bigskip
\noindent{\bf NB :} the case $t=0$ is included in $\v\,\frac{z}{t} \leq 0$.
\bigskip

We can now finish the proof of Lemma \ref{FV}. 
The same computation, with $v_2$ and $b$ instead of $v_1$ and $a$, gives the values of $\ga^{-i}v_2$ for any $i \in \N$. 
It is then easy to compute $F(k,k')$ as given in Lemma \ref{FV}, for $k$ and $k'$ in K, and to check formula (\ref{FK}).

\subsection{Some test-vectors}

On the one hand, from the expression of $F$ given by Lemma \ref{FV} we deduce 
$$\Bigl( \Psi(F) \Bigr)(v_3) = A \cdot \ell(v_1'\otimes v_2'\otimes v_3)$$
On the other hand, from the relation $F=\ext(f)$ we deduce 
$$\Bigl( \Psi(F) \Bigr)(v_3) = \Bigl( \Phi(f) \Bigr)(v_3) $$
and from equation \ref{phi-f-v_3} in Section \ref{phi-f} we get 
$$\Bigl( \Psi(F) \Bigr)(v_3)  \not= 0 $$
Then 
$$\ell(v_1'\otimes v_2'\otimes v_3) \not = 0$$
and $v_1'\otimes v_2'\otimes v_3$ is a test vector for $\ell$. 
We are going to simplify it. We can deduce from lemma \ref{FV} that 
$$ F =   A \cdot
\left\{ ab \cdot \ga^{-(n-1)}v_1 \otimes\gb v_2  - a \cdot \ga^{-(n-1)}v_1 \otimes v_2 
- b \cdot \ga^{-n}v_1 \otimes\gb v_2  + \ga^{-n}v_1 \otimes v_2 \right\} .$$ 

\bigskip
If $n\geq 2$, we write 
$$\Bigl( \Psi(F) \Bigr)(v_3) = 
A \cdot \bigl( ab \cdot \langle \gb \psi_{n-2}, v_3 \rangle -a \cdot \langle \psi_{n-1}, v_3\rangle-
b \cdot\langle \gb\psi_{n-1}, v_3\rangle + \langle \psi_n, v_3\rangle \bigr)$$ 
where, for $m$ in $\{n-1,n-2,n\}$ 
$$\psi_m \quad \left\{ \begin{matrix}
V_3 & \longrightarrow & \C \hfill \cr
v & \longmapsto & \ell(\ga^{-m}v_1\otimes  v_2 \otimes v) \cr
\end{matrix}  \right.$$
Since $\ell$ is $G$ invariant, $\psi_m$ is an element of $\widetilde{V_3}$ which is invariant by the action of
$$R_m \cap K =I_m $$
But $\widetilde{\pi_3}$ has conductor $n$ so  
$$\psi_{n-2}=\psi_{n-1} = 0$$
and 
$$\Bigl( \Psi(F) \Bigr)(v_3) = A \cdot  \psi_3( v_3) = A \cdot \ell(\ga^{-n}v_1\otimes  v_2 \otimes v_3)$$ 
then $$\ell(\ga^{-n}v_1\otimes  v_2 \otimes v_3) \not=0 $$
and $\ga^{-n}v_1\otimes  v_2 \otimes v_3$ is a test vector for $\ell$. 

\bigskip
If $n = 1$, only the two terms in the middle vanish and we get
$$\Bigl( \Psi(F) \Bigr)(v_3) = A \cdot \bigl( ab \cdot \ell( v_1\otimes \gb v_2 \otimes v_3) + \ell(\gb v_1\otimes v_2 \otimes v_3)\bigr)$$ 
Now, take
$$g=\begin{pmatrix}0 & 1 \cr \PF & 0\cr\end{pmatrix} .$$
On the one hand 
$$g \gb = \begin{pmatrix} 0 & 1 \cr 1 & 0 \cr \end{pmatrix}.$$
and this matrix is in $K$ so 
$$g \gb v_1 = v_1$$ 
On the other hand
$$\begin{pmatrix} 0 & 1 \cr \PF & 0 \cr \end{pmatrix}= 
\begin{pmatrix} \PF & 0 \cr 0 & \PF  \cr \end{pmatrix}
\begin{pmatrix} \PF^{-1}  & 0 \cr 0 & 1 \cr \end{pmatrix}
\begin{pmatrix} 0 & 1 \cr 1 & 0 \cr \end{pmatrix}.$$
The first matrix belongs to the center of $G$, the second one is precisely $\gb$ and 
the third one is in $K$, so 
$$g v_2 =  \gb v_2.$$
Then
\begin{align*}
\Bigl( \Psi(F) \Bigr)(v_3)   
 = & A \cdot \bigl\{ ab \cdot \ell( g \gb v_1\otimes g v_2 \otimes v_3) + \ell(\gb v_1\otimes v_2 \otimes v_3)\bigr\}\hfill \\
 = & A \cdot \bigl\{ ab \cdot \ell( \gb v_1\otimes  v_2 \otimes g^{-1}v_3) + \ell(\gb v_1\otimes v_2 \otimes v_3)\bigr\}\hfill\\
 = & A \cdot \ell\bigl( \gb v_1\otimes v_2 \otimes v_3' \bigr)\hfill\\
\end{align*}
with 
$$ v_3' = ab \cdot (g^{-1}v_3) + v_3$$
The linear form 
$$\psi_1 \quad \left\{ \begin{matrix}
V_3 & \longrightarrow & \C \hfill \cr
v & \longmapsto & \ell(\gb v_1\otimes  v_2 \otimes v) \cr
\end{matrix}  \right.$$
is a non zero element of ${\widetilde{V_3}}^{I_1}$, so it is a new vector in $\widetilde{V_3}$.
It is known from \cite{B}  that a new vector in $\widetilde{V_3}$ does not vanish on on $V_3^{I_1}$ that is, it does not vanish on a new vector of $V_3$ : 
$$\ell(\gb v_1 \otimes v_2 \otimes v_3) \not= 0 $$
and $\gb v_1 \otimes v_2 \otimes v_3$ is a test vector for $\ell$. 
This concludes the proof of Theorem \ref{vt-00n}.

\vskip1cm
{\bf NB :} it is easy to deduce from Theorem \ref{vt-00n}  that $ v_1 \otimes \gbn v_2 \otimes v_3$ also is a test vector for $\ell$. 
Take $$g=\begin{pmatrix}0 & 1 \cr \PF^n & 0\cr\end{pmatrix}.$$
Then
$$ g \gbn  v_1 = v_1 \quad ,\quad g v_2 = \gbn v_2$$
and
$$\ell(v_1 \otimes \gbn v_2 \otimes g v_3) = \ell( g \gbn v_1 \otimes g v_2 \otimes g v_3) = \ell(  \gbn v_1 \otimes  v_2 \otimes  v_3) \not=0.$$
So the linear form 
$$\left\{ \begin{matrix}
V_3 & \longrightarrow & \C \hfill \cr
v & \longmapsto & \ell(v_1\otimes  \gbn v_2 \otimes v) \cr
\end{matrix}  \right.$$
is not zero. Being $\In$-invariant, it is a new vector in $\widetilde{V_3}$, which does not vanish on $v_3$ : 
$$\ell(v_1 \otimes \gbn v_2 \otimes  v_3)\not=0.$$

\end{document}